\documentclass[11pt]{article}
\usepackage[T1]{fontenc}
\usepackage[utf8]{inputenc}

\usepackage{lmodern,amsmath,amssymb,amsthm,mathtools}
\usepackage[margin=1in]{geometry}
\usepackage{microtype,booktabs,array,enumitem,listings,needspace}
\usepackage{longtable}
\IfFileExists{bbm.sty}{\usepackage{bbm}}{}
\usepackage{todonotes} 
\usepackage[colorlinks=true,linkcolor=blue,citecolor=blue,urlcolor=blue]{hyperref}
\newtheorem{theorem}{Theorem}[section]
\newtheorem{lemma}[theorem]{Lemma}
\newtheorem{proposition}[theorem]{Proposition}
\newtheorem{corollary}[theorem]{Corollary}
\theoremstyle{definition}
\theoremstyle{remark}\newtheorem{remark}[theorem]{Remark}
\newcommand{\NIM}{\operatorname{NIM}}
\newcommand{\Pset}{\mathcal P}
\newcommand{\Nset}{\mathcal N}
\newcommand{\Zge}{\mathbb Z_{\ge0}}
\newcommand{\BA}{\operatorname{BA}}

\newcommand{\val}{\operatorname{bin}}

\ifdefined\mathbbm
\def\uu{\mathbbm{1}}
\else
\def\uu{\mathbf{1}}
\fi
\let\es\varnothing
\let\ld\lambda
\let\sm\setminus

\def\A{\mathcal{A}}

\def\L{\mathsf{L}}

\def\NC{\mathsf{NC}}
\def\FNC{\mathsf{FNC}}

\def\tdm{\widetilde m}
\def\tS{\widetilde S}
\def\poly{\mathop{\mathrm{poly}}\nolimits}
\def\sdif{\mathbin{\triangle}}
\def\TR{T_{\text{rank}}}
\def\TS{T_{\text{Smith}}}
\allowdisplaybreaks
\setlist[enumerate]{itemsep=.35em,topsep=.5em}
\newcommand{\NN}{\mathbb Z_{\geq0}}

\hypersetup{pdftitle={Computing the Remoteness Function of Exact NIM(n,k) with n <= 3k},
pdfauthor={Vladimir Gurvich, Mariya Naumova, Artem Parfenov, Mikhail Vyalyi}}
\title{Computing the Remoteness Function of Exact $\NIM(n,k)$ with $n \leq 3k$}
\author{
Vladimir Gurvich\thanks{National Research University Higher School of Economics and RUTCOR, Rutgers University, New Jersey, USA, \texttt{vladimir.gurvich@gmail.com}.}\and
Mariya Naumova\thanks{Rutgers Business School, Rutgers University, Piscataway, NJ, United States, \texttt{mnaumova@business.rutgers.edu}.}\and
Artem Parfenov\thanks{National Research University Higher School of Economics, Moscow, Russia, \texttt{dunno\_o@icloud.com}.} 
\and
Mikhail Vyalyi\thanks{Federal Research Center “Computer Science and Control” of the Russian Academy of Sciences; National Research University Higher School of Economics; Moscow Institute of Physics and Technology;  \texttt{vyalyi@gmail.com}.}}
\date{\today}

\begin{document}
\maketitle
\begin{abstract}
We obtain an explicit formula for the Smith remoteness function 
of the impartial game Exact $\NIM(n,k)$ when $2k\leq n \leq 3k$  
under the normal-play convention (the case $n<2k$ is trivial). 
A legal move in this game decreases exactly $k$ of $n$ piles 
of tokens by any positive amounts.  
Even characterizing the P-positions of $\NIM(5,2)$ 
was an open problem for at least ten years.
\newline 
For each position we consider the corresponding 
Bouton matrix, introduce the concept of its rank,  
and define capacity as the maximal possible number of piles that 
can be decreased while preserving the rank.
Then, P-positions are exactly those of capacity smaller than $k$.
Generalizing this result, we also express the remoteness 
function in the same terms.  
\newline
We provide an algorithm computing 
the rank and capacity in polynomial time 
even if $n$ and $k$ are part of the input. 
For fixed $k$, we provide a finite-state automaton recognizing 
the P-positions and nondeterministic finite-state transducers 
that output the rank and an optimal move 
(the winner wins as quickly as possible; 
the loser resists as long as possible). 
From this we conclude that, for Exact $\NIM(n,k)$ with $k=O(1)$ and $n\leq 3k$, P-membership belongs to logspace uniform $\NC^1$, while computing the remoteness function and an optimal move belongs to $\FNC^1$; these problems can therefore be solved in logarithmic space.
\end{abstract}

\section{Introduction and background}\label{sec:background}

\subsection{Outcomes, Sprague--Grundy and Smith remoteness functions}
In a finite impartial game under the normal-play convention, a P-position is losing for the
player about to move, whereas an N-position is winning  
(in other words, the Previous player and the Next player, respectively, win).  
A set of positions is the P-set precisely when 
it contains every terminal, no move connects
two of its members, and every position outside it has a move into it.

There are two distinct numerical refinements defined recursively. 
The Sprague--Grundy (SG) function 
(also called the nim-value) introduced in~\cite{Sprague35,Sprague37,Grundy}
is the mex of the values of the options and vanishes exactly on P-positions.
The Smith remoteness function $R(x)$  
(also called the remoteness value or simply the remoteness, for short) instead records 
the length of optimal play when the winner wins as quickly as possible, and the loser
resists as long as possible \cite{S66}. Its even values occur on P-positions,
and its odd values on N-positions. The precise recurrence is given below;  
see Eq.~\eqref{eq:remoteness}. 
An \emph{optimal move} is a legal move $x\to y$ satisfying $R(y)=R(x)-1$.

These two functions play analogous roles for two different ways of combining games.
The first one is related to so-called  
{\em disjunctive compounds} of games. 
The SG-value for a disjunctive compound is the NIM-sum 
(bitwise XOR) of the SG-values of the components.
These results were obtained in 1901 by Bouton 
\cite{B02}  for the special case of NIM and 
then extended to arbitrary impartial games 
by Sprague \cite{Sprague35,Sprague37} and Grundy \cite{Grundy}.\footnote{Bouton was also the first to consider 
both the mis\`ere and normal conventions, for the special case of Nim.  
He was motivated to do so
because Harvard students mostly played the mis\`ere version.}

The remoteness function $R(x)$ is related to
\emph{conjunctive compounds}. In the convention considered here,
a legal move in a conjunctive compound consists of making one
legal move in every component. Consequently, the compound is
terminal as soon as at least one component is terminal.
Under normal play, the player whose turn it is at such a position
loses. If $x^{(1)},\ldots,x^{(m)}$ are positions in the component
games and $R_1,\ldots,R_m$ are their Smith remoteness functions,
then the remoteness of their conjunctive compound is
\[
  R_{\wedge}\bigl(x^{(1)},\ldots,x^{(m)}\bigr)
  =\min_{1\leq i\leq m}R_i\bigl(x^{(i)}\bigr).
\]
This result was obtained in 1966 by Smith~\cite{S66};
see definitions and details in~\cite[Section~1.4]{BGMV24}.

\subsection{Related work}

The binary description of losing positions in Nim, established 
in 1901 by Bouton~\cite{B02}, is one of the earliest examples of a complete
mathematical solution of an impartial game.  
Bouton's classical solution identifies P-positions 
by the vanishing of the bitwise XOR of the pile sizes \cite{B02}. 
Equivalently, every column of the binary pile matrix has even weight. 
This criterion already suggests that 
the appropriate representation of a position is a matrix of digits 
rather than the magnitudes of the piles.

In 1910, Moore~\cite{M10} generalized Nim 
by allowing players to decrease between 1 and $k$ piles in one move 
rather than exactly 1. 
He showed that P-positions\footnote{which have the SG-value 0.} 
of this game, $\NIM_M(n,k)$, are exactly those 
whose binary Bouton matrix has column weights divisible by $k+1$.
In 1980, Jenkyns and Mayberry \cite{JenkynsMayberry}  
obtained a similar formula for positions of SG-value 1. 
In the same paper they also proved an explicit formula 
for the SG function of $\NIM_M(n,k)$ with $n = k+1$. 
Already for $\NIM_M(4,2)$ and the SG-value 2,
no polynomial-time criterion for recognizing these positions is known. 

Recently, Boros, Gurvich, Makino, and Vyalyi~\cite{BGMV24} 
studied the remoteness function of several well-known impartial games.
For Moore Nim they gave a simple explicit formula
for the remoteness of P-positions (which have even remoteness), 
as well as for optimal moves from these positions.
In contrast, computing Smith-optimal moves 
from N-positions to P-positions is NP-hard 
when $n$ and $k$ are part of the input.

Here we study  \emph{Exact} $\NIM(n,k)$. 
A legal move in this game must
strictly decrease  exactly $k$ piles. Different selected
piles may lose different positive numbers of tokens.
Gurvich, Parfenov, and Vyalyi~\cite{GPV23} investigated relations
between the P-positions of Exact $\NIM(5,2)$ and those of
$\NIM_M(4,2)$, for which Moore's criterion applies. Their
experimental study provided evidence of these relations but did
not obtain a complete characterization of the P-positions of
Exact $\NIM(5,2)$.

Boros, Gurvich, Ho, Makino, and Mur\v{s}i\v{c} determined the
SG function of Exact $\NIM(n,k)$ for $n\leq2k$ in~\cite{BGHMM18};
the broader hypergraph framework was developed in~\cite{BGHMM19}.
When $n<2k$, the SG-value of a position equals its \emph{height},
that is, the maximum length of a play from that position to a
terminal. Consequently, every legal move strictly decreases the
SG-value. Games with this property are called \emph{SG-decreasing};
see~\cite{BGHMM24} for their structural theory and hypergraph
compounds. The boundary case $n=2k$ with $k\geq2$ has a different
structure: its SG function is given by the JM formula, named in
honor of Jenkyns and Mayberry~\cite{JenkynsMayberry}; see
\cite{BGHMM18,BGHMM19} for the formula and its applications.
The exception $k=1$ is ordinary two-pile Nim, whose SG function is
given by Bouton's XOR formula. These results do not provide a
general SG formula for $n>2k>2$. 

In the hypergraph interpretation developed in~\cite{BGHMM19},
the vertices represent piles, and a move selects a hyperedge and
strictly decreases all its coordinates. Exact $\NIM(n,k)$
corresponds to the hypergraph of all $k$-subsets of the set of
$n$ piles, whereas Moore's $\NIM_M(n,k)$ corresponds to all
nonempty subsets of size at most $k$. This framework places the
boundary case $n=2k$ within a broader structural theory and
explains why the same JM formula applies to several classes of
hypergraphs. Structural information does not necessarily make
value computation easy: computing the height and the SG-value
is NP-hard even for intersecting hypergraphs, whose associated
Nim games are SG-decreasing, when the hypergraph is part of the
input; see~\cite[Theorem~9 and Corollary~2]{BGHMM24}.
This evaluation problem is distinct from recognizing whether a
given hypergraph is SG-decreasing.

The remoteness function for the boundary case $n=2k$ is also covered
by the results of Boros, Gurvich, Makino, and Vyalyi~\cite[Theorem 11]{BGMV24}
for minimally transversal-free hypergraphs. Theorem~\ref{thm:matrix}
recovers this case and extends the outcome and remoteness characterization
to $2k<n\leq3k$ for $k\geq2$.

\medskip

Requiring each selected pile to lose exactly one token defines
Exact \emph{Slow} $\NIM^1(n,k)$, a different game: every move
decreases exactly $k$ nonempty piles by one token each.
The slow variants were introduced by Gurvich and Ho
in~\cite{GH15}. Subsequently, the SG-values of
$\NIM^1(3,2)$ and $\NIM^1(4,2)$ were computed for both the
normal and mis\`ere versions in~\cite{GHHC20}, and the
normal-play P-positions of $\NIM^1(5,2)$ and $\NIM^1(6,2)$
were characterized in~\cite{CGKPV21}. More recently, an explicit
formula for the normal-play SG function of $\NIM^1(n,k)$
with $n=k+1$ and $k\geq2$ was obtained in~\cite{GN26}.

The normal-play remoteness function for the family $n=k+1$
was determined in~\cite{GMMV24}. Notably, an optimal move from
every nonterminal position can be chosen by the M-rule, without
first deciding whether the position is P or N. Thus, the move
rule is simpler to apply than a separate characterization of
the P-positions. The M-rule, its generalizations, and the
associated discrete dynamics were studied in~\cite{GN24}.

In~\cite{GMMN23}, the M-rule was investigated for the mis\`ere
version of $\NIM^1(n,k)$ with $n=k+1$. The results are partial:
computations show that the rule is optimal at many positions
and identify structured exceptional families. Formulas were
proposed for the exceptions found computationally, but their
exhaustiveness was not proved.

In~\cite{GMN26}, the generalized M-rule, or GM-rule, was
investigated for normal-play $\NIM^1(n,k)$ with $k=2$.
For $\NIM^1(3,2)$, the M-rule and GM-rule coincide and give
an optimal Smith move from every nonterminal position.
For $\NIM^1(4,2)$ and $\NIM^1(5,2)$, computations reveal
sparse, structured exceptions to the GM-rule, and simple
formulas are proposed to describe all known exceptions.
These descriptions are computationally supported rather than
proved exhaustive classifications.

Generalizations of Exact Slow Nim, including P-position
characterizations for infinite families of Slow $A$-Nim,
were studied in~\cite{DH26}. These results concern a restricted
move set and do not, by themselves, establish the optimality
of slow moves in unrestricted Exact Nim.

\subsection{Our contribution} 

Let us note that the case $n < 2k$ is trivial. 
If  $n < k$ then every position is terminal. 
If $k \leq n < 2k$ then from every nonterminal 
position a terminal one can be reached in one move: 
just eliminate all tokens from any $k$ nonempty piles.  
Thus, in this case the remoteness is 1 for each 
nonterminal position, and it is 0 for every terminal one, by definition. 

We compute the remoteness function of Exact $\NIM(n,k)$  for $2k \leq n \leq3k$. 
To do so, we introduce the {\em rank} and {\em capacity} of a position. 
The latter is defined as the maximal possible number of piles 
that can be decreased while keeping the rank unchanged.
Then, P-positions are exactly those whose capacity is smaller than $k$, and 
the remoteness function is also expressed in the same terms; see Theorem~\ref{thm:matrix} below. 

Even the P-positions of Exact $\NIM(5,2)$ were not known. 
This was mentioned as an important open problem in several papers~\cite{BGHMM18,GPV23}.

We also study the algorithmic complexity of basic problems related to $\NIM(n,k)$: 
determining the outcome, computing the value of the remoteness function, and finding optimal moves. 
The crucial result is a greedy algorithm for computing the capacity; 
see Lemma~\ref{lem:capacity-count} below. 
This algorithm can be implemented in polynomial time. 
It places all basic algorithmic problems 
for Exact $\NIM(n,k)$ with $2k \leq n \leq 3k$ 
in polynomial-time complexity classes; see Theorem~\ref{thm:nk-input} below.\footnote{In contrast, for Moore Nim with $n$ and $k$ as part
of the input, computing a Smith-optimal move from an N-position
is NP-hard. More precisely, given an N-position $x$ and an
integer bound $B_0$, deciding whether $R_M(x)\leq B_0$ is
NP-complete, where $R_M$ denotes the remoteness function of
Moore Nim; see~\cite{BGMV24}.}

For a fixed $k=O(1)$, we provide stronger results. 
Based on Lemma~\ref{lem:capacity-count}, we provide 
 a finite-state automaton recognizing the set of P-positions and nondeterministic finite-state transducers that output the rank and an optimal move; see Propositions~\ref{prop:scan} and~\ref{Smith-optimal:NFST}. 
 From this we derive that basic algorithmic problems 
 related to Exact $\NIM(n,k)$  with $k=O(1)$, $2k\leq n\leq 3k$, can be solved in linear time (Theorem~\ref{thm:Smith-move-linear-time-linear-space}); P-membership belongs to uniform $\NC^1$, while computing the remoteness function and an optimal move belongs to $\FNC^1$ (Theorem~\ref{thm:fixed-circuits}). These problems can therefore be solved in logarithmic space.

\section{Main results}
\subsection{Preliminaries} 
Fix $n\geq k\geq1$ such that $2k\leq n \leq 3k $. Throughout the paper
we use the notation $r=n-2k$ (so $0\leq r\leq k$).

A position of Exact $\NIM(n,k)$ is a vector $x=(x_1,\ldots,x_n) \in\Zge^n$.
A legal  move $x\to y$ strictly decreases exactly $k$ entries and 
leaves the others unchanged. Reductions can be
unequal and larger than one. A terminal has fewer than $k$ positive piles. 

A P-position is a position from which the player
about to move loses against correct play, and an N-position is a
position from which that player can force a win.
We write $\Pset$ and $\Nset$ for these sets.
Since the sum of the coordinates strictly decreases along every move, 
backward induction on this sum shows that 
the sets of P-positions and N-positions  are characterized by the fact that no
legal move joins two P-positions and every N-position has a move to a 
P-position, with every terminal position belonging to $\Pset$.
The remoteness function is defined by
\begin{equation}\label{eq:remoteness}
R(x)=\left\{\begin{aligned}
&0,&&x\text{ is terminal},\\[1mm]
&1+\min\{R(y):x\to y,\ y\in\Pset\},&&x\in\Nset,\\[1mm]
&1+\max\{R(y):x\to y\},&&x\in\Pset\text{ is nonterminal}.
\end{aligned}\right.
\end{equation}
Backward induction shows that $R$ is even on $\Pset$ and odd on $\Nset$.
We call a move \emph{optimal} when it decreases $R$ by one, so that
the term refers both to a fastest forced win and to a longest
possible resistance, according to the outcome of the initial position.


\subsection{The Bouton matrix}

Choose $L\geq1$ with every $x_i<2^L$. 
The Bouton matrix of a position $x$ is a matrix $B(x)$ of size $n\times L$ such that 
\[
 B(x)=(X_{ij})_{\substack{1\leq i\leq n\\1\leq j\leq L}},
 \qquad x_i=\sum_{j=1}^L X_{ij}2^{L-j}.
\]
Thus the most significant column is on the \emph{left}. Leading zero
columns are allowed. For a binary word $u=(u_1,\ldots,u_L)$, put
$\val(u)=\sum_j u_j2^{L-j}$. Lexicographic order and numerical order agree
for words of this common length. For example,
\[
 B(3,5,8)=\begin{pmatrix}
 0&0&1&1\\0&1&0&1\\1&0&0&0
 \end{pmatrix}.
\]
In a strict decrease of a row, the leftmost change is
$1\to0$; entries to the right may then be arbitrary.

In what follows we identify Bouton matrices and positions. In particular, we use the notations $x$ and $X=B(x)$ interchangeably for positions. Also, we identify rows and columns with their indices, for brevity.

\subsection{Rank and capacity} 

We use bitwise Boolean operations with integers. Given a representation of integers in binary, the result of such an operation is obtained by applying it to the corresponding bits.

Let $(x_1,x_2,\dots, x_n)$ be a position.
Define the \emph{rank} of $x$ by
\begin{equation}\label{eq:rank}
b=b(x)=
\min_{\substack{I\subseteq\{1,\ldots,n\}\\|I|=r+1}}
\bigvee_{i\in I}x_i,
\end{equation}
where the minimum is taken in ordinary numerical order. Recall that $r= n-2k$.

\begin{proposition}\label{b-terminal}
  If $x$ is terminal then $b(x)=0$.
\end{proposition}
\begin{proof}
  In a terminal position there are fewer than $k$ nonzero entries. So, there are at least $2k+r - k+1 = k+r+1$ zero  rows. Bitwise OR of $r+1$ of them is zero. 
\end{proof}

For a pair of positions $x$, $y$ we say that  $x\xrightarrow{d} y$ is a $d$-pseudomove if 
$  y_i\leq x_i$ for all $1\leq i\leq n$ and   $|\{i:y_i<x_i\}|=d$. A $0$-pseudomove $x\xrightarrow{0} x$ is allowed. Thus,  $0\leq d\leq n$.
If $b(y)= b(x)$ then we say that $x\xrightarrow{d} y$ is a \emph{conservative} $d$-pseudomove.  

The \emph{capacity} $A(x)$ of a position $x$ 
is the maximal $d$ such that there exists a conservative $d$-pseudomove from $x$.

\subsection{Main theorem}

We show that $A(x)$ determines 
the outcome of the game and its remoteness function. 

\begin{theorem} 
\label{thm:matrix}
For Exact $\NIM(n,k)$ with $2k\leq n\leq 3k$, 
a position $x$ is a P-position if and only if $A(x)<k$.
Moreover, let $\varepsilon=\uu[A(x)\geq k]$ be the indicator function of N-positions. Then, 
\begin{equation}\label{eq:matrixR}
 R(x)=2b(x)+\varepsilon.
\end{equation}
\end{theorem}

We start with a sketch of the proof. 
Let $P_A = \{x: A(x)<k\}$ and $N_A = \{x: A(x)\geq k\}$. Note that for a terminal position $A(x)\leq|\{i:x_i>0\}|<k$. So, all terminal positions belong to $P_A$. 
To prove the first claim of Theorem~\ref{thm:matrix}, 
we need to show that there are no moves from $P_A$ to $P_A$ and 
for each $x\in N_A$ there exists a move $x\to y$ such that $y\in P_A$. 
To prove the second claim, we need more. 
Namely, we are going to show that
\begin{enumerate}[label=\textup{(\roman*)}]
\item\label{nonincrease} the rank is strictly decreasing 
on a move from a position in $P_A$ (Lemma~\ref{lem:rank} below); 
\item\label{nondecrease} the rank is non-decreasing on a move 
to a position in $P_A$ (Lemma~\ref{lem:rank} below); 
\item\label{Pmove} for each nonterminal $x\in P_A$ there exists a move $x\to y$ such that $b(y) = b(x)-1$ (Lemma~\ref{lem:P-Smith} below);
\item\label{Nmove} for each $x\in N_A$ there exists a move $x\to y$ such that $y\in P_A$ and $b(x) = b(y)$ (Lemma~\ref{NPrank-preserve} below).
\end{enumerate}

Items~\ref{nonincrease} and~\ref{nondecrease} show that  there is no move from $P_A$ to $P_A$.

Item~\ref{Nmove} provides a move from any position in $N_A$ to a position in $P_A$. 
It completes the proof of the first part of Theorem~\ref{thm:matrix}. Thus, $P_A=\Pset$, $N_A=\Nset$.
For the second part, note that the items above imply a recurrence on the rank:   
\begin{equation}\label{eq:rank-rec}
b(x)=\left\{\begin{aligned}
&\min\{b(y):x\to y,\ y\in P_A\},&&x\in N_A,\\[1mm]
&1+\max\{b(y):x\to y\},&&x\in P_A\ \text{ is nonterminal}.
\end{aligned}\right.
\end{equation}
By induction on $\sum_i x_i$, we verify the recurrence~\eqref{eq:remoteness} for the function $R(x) = 2b(x)+\varepsilon$ using the recurrence~\eqref{eq:rank-rec}. For a terminal position $x$, we have  $R(x)=0$, $b(x)=0$ (by Proposition~\ref{b-terminal}), and $A(x)<k$.
If $x\in\Nset=N_A$ then $R(x) = 1+\min\{2b(y): x\to y,\ y\in\Pset\}=2b(x)+1$  by the first line of~\eqref{eq:rank-rec}.
If $x\in\Pset$ and $x$ is nonterminal then  $R(x)=1+\max\big\{2b(y)+1: x\to y\big\}=2b(x)$ by  the second line of~\eqref{eq:rank-rec}.

In the next subsections we implement this plan step by step. First, however, we introduce several important concepts related to the rank.

\subsection{The scan for the rank word}

The minimum in~\eqref{eq:rank} can be found by a greedy procedure  scanning through the columns of $X$ from left (the most significant) to right (the least significant).
The scan generates sequences $I(j)$, $Z(j)$, $E(j)$, where $1\leq j\leq L$,  of subsets of the row set $\{1,\ldots,n\}$.

Initially, the \emph{active row set} is $I(0)=\{1,\ldots,n\}$ and the \emph{marked row set} is $E(0)=\es$.

 At column $j$, let $Z(j)=\{i\in I(j-1):X_{ij}=0\}$. Then, for $1\leq j\leq L$, set
\begin{equation}\label{eq:word-scan}
\left\{ \begin{aligned}
 \lambda_j&=0, &&I(j) = Z(j), && E(j) =\{1,2,\ldots, n\}\sm Z(j),&&\text{if } |Z(j)|\geq r+1,\\
 \lambda_j&=1,&&I(j) = I(j-1), && E(j) = E(j-1),&&\text{if } |Z(j)|\leq r.
 \end{aligned}\right.
\end{equation}
In words, we set $\ld_j=1$ if and only if there are at most $r$ zeros in active  rows of the current column.  
If $\ld_j=0$, the next set of active rows is restricted to the rows
that have digit 0 in the current column. The other rows become marked.
Finally, we set $I = I(L)$ and $E= E(L)$.

\begin{proposition}\label{E-bound}
    $I(j) \cap E(j) =\es$, $I(j)\cup E(j) = \{1,\ldots, n\}$, and $|I|\geq r+1$, $|E|<2k$.
\end{proposition}
\begin{proof}
    It follows immediately from the definition.
\end{proof}

The word  $\ld(X) $ is called the \emph{rank word}.
We say that column $j$ is \emph{active} if $\ld_j=1$.

\begin{lemma}\label{lem:word}
For a position $(x_1,\dots, x_n)$ and its Bouton matrix, $b(X) =\val \big(\ld(X)\big) $.
Every $r+1$ active rows attain the minimum in~\eqref{eq:rank}.
\end{lemma}
\begin{proof}
By induction on $j$, we prove both  claims for all prefixes of the Bouton matrix.
Initially every set of cardinality $r+1$  is possible. If $\ld_j=1$ then the bitwise OR of any set of $r+1$ entries in active rows and column $j$ is 1.    If $\ld_j=0$ then the bitwise OR of any set of $r+1$ entries in $Z(j)=I(j)$ is 0.
\end{proof}

In the example $(3,5,8)$, for $k=1$ (so $r=1$, too), we get $\lambda=0111$ and $b=7$.

\subsection{Normalization   and the normalized matrix}

The \emph{normalized} matrix $C=C(X)$ is constructed as follows.  In an active row $i$, it coincides with $X$.
In a marked row, find the leftmost column $\ell(i)$
where $X_{i \ell(i)}=1$ but $\lambda_{\ell(i)}=0$. Note that  before processing the column $\ell(i)$ the row $i$ is active and after processing the column $\ell(i)$ the row $i$ is marked. So, this column is called the \emph{marking column} for the row $i$. Leave  entries of the row to the
left of $\ell(i)$ unchanged, and replace its entire suffix starting at $\ell(i)$
by the corresponding suffix of $\lambda$:
\begin{equation}\label{eq:normal}
 C_{ij}=\left\{\begin{aligned}
 X_{ij},&&&j<\ell(i),\\
 \lambda_j,&&&j\geq \ell(i).
 \end{aligned}\right.
\end{equation}
For brevity, we will use the notation $\ell(i)=\infty$ for active rows.

\begin{lemma}\label{lem:normal}
The normalized matrix has the following properties:
\begin{enumerate}[label=\textup{(\roman*)}]
\item\label{i} If $\lambda_j=0$ then column $j$ of $C$ is zero.
\item\label{ii} If $\lambda_j=1$ then column $j$ of $C$ has weight at least $n-r$.
\item\label{iii} For every set of $r+1$ rows of $C$, their  bitwise OR is equal to  $\lambda$.
\end{enumerate}
\end{lemma}
\begin{proof}
  If $\ld_j=1$ then active rows in $I(j)$ have at most $r$ zeros in column $j$. The marked rows in $E(j)$ have digit 1 by definition of the normalized matrix.
   If $\lambda_j=0$, then active rows in $I(j)$ have 0 in column $j$ by definition~\eqref{eq:word-scan}. 
   The marked rows in $E(j)$ have 0 by definition of the normalized matrix.
   This proves (i) and (ii). They imply (iii), since $\bigvee_{i\in I}C_{ij}=0 $, $|I|=r+1$, implies  $\ld_j=0$ and  $\bigvee_{i\in I}C_{ij}=1 $, $|I|=r+1$,  implies  $\ld_j=1$.
\end{proof}

\begin{proposition}\label{AgeqE}
  $A(X)\geq |E|$.
\end{proposition}
\begin{proof}
Each marked row strictly decreases lexicographically after normalization. Also, normalization does not change the rank. This implies the claim.
\end{proof}

\begin{corollary}\label{many-zeroes}
  If  $X\in P_A$ then $X$ has at least $k+r+1$ active rows.
\end{corollary}
\begin{proof}
  From Proposition~\ref{AgeqE}, $|E|\leq A(X)<k$. Therefore   $|I|= 2k+r-|E|> k+r $.
\end{proof}





\subsection{Rank changes after a move}

\begin{lemma}\label{lem:rank}
Let  $X\to Y$ be a move in $\NIM(n,k)$, $2k\leq n\leq 3k$. If $X\in P_A$ then $b(Y)< b(X)$. If $Y\in P_A$ then $b(Y)\geq b(X)$.    
\end{lemma}
\begin{proof}
   Active rows coincide in a Bouton matrix and in its normalization.
   For a position from $P_A$, there are at least $k+r+1$ active rows due to Corollary~\ref{many-zeroes}. Since a move changes exactly $k$ rows, at least $r+1$ active rows remain unchanged. Bitwise OR of $r+1$ unchanged rows equals $\ld (X)$ by Lemma~\ref{lem:normal}. So, $b(X)=\val\big(\ld(X)\big)\geq b(Y)$. The equality is impossible by the definition of capacity.

  Suppose $Y\in P_A$. Again,   $\ld(Y)$ is the bitwise OR of $r+1$   unchanged  rows that are active in $Y$. Thus, $b(X)\leq b(Y)$.
\end{proof}

\subsection{An optimal move from a P-position}

\begin{lemma}\label{lem:P-Smith}
Let $x\in P_A$ be a nonterminal position. Then there exists  a move $x\to y$ such that $b(y)=b(x)-1$.   
\end{lemma}
\begin{proof}
  By Corollary~\ref{many-zeroes}, there are at least $k+r+1$ active rows of the Bouton matrix $X = B(x)$. Active rows of $X$ and $C(X)$ coincide. By Lemma~\ref{lem:normal}, the bitwise OR of any set of $r+1$ active rows equals $\ld(X)$. Thus, $b(x)=0$ implies that at least $k+r+1$ rows are zero and the position is terminal.

  Now assume that $b(x)>0$. Let $s$ be the length of the maximal zero suffix of the rank word $\ld(X)$.
  By Lemma~\ref{lem:normal}, active rows also have zero suffixes of length $s$, and at most $r$ of them have zero in column $L-s$, that is, column $s+1$ counting from the right. Hence at least $k+1$ active rows have digit 1 in that column. Select $k$ of them and make a slow move $x\to y$ that decreases each selected pile (row) by 1.

  The suffix of each selected row changes as follows:
\[
u\,1\,\underbrace{00\cdots0}_{s}
       \quad\longmapsto\quad
       u\,0\,\underbrace{11\cdots1}_{s}.
\]
Consider any set of $r+1$ rows. If all these rows are unchanged, their bitwise OR is unchanged and its value is at least $b(x)$ by~\eqref{eq:rank}. Otherwise, their OR has the same prefix strictly to the left of column $L-s$ as before the move, and its last $s$ entries are ones. Since the original OR has value at least $b(x)$, this unchanged prefix is lexicographically at least the corresponding prefix of $\ld(X)$. If it is larger, the new OR is greater than $b(x)-1$; if it is equal, the smallest possible suffix is $0\,\underbrace{11\cdots1}_{s}$, which is exactly the suffix of the binary representation of $b(x)-1$. Thus, every such OR has value at least $b(x)-1$, and $b(y)\geq b(x)-1$. But $b(y)<b(x)$ by Lemma~\ref{lem:rank}. Therefore $b(y)=b(x)-1$.
\end{proof}

\subsection{An optimal move from an N-position}
\label{sec:simple-N}

Lemma~\ref{lem:rank} implies that any move $x\to y$ to a position  $y\in P_A$ does not decrease the rank. So, $b(x)$ is the smallest possible value of $b(y)$. 

\begin{lemma}\label{NPrank-preserve}
  Let $X\in N_A$. Then there exists  a move $X\to Y$, $Y\in P_A$ that  preserves the rank: $b(X)=b(Y)$. 
\end{lemma}

Before the proof of Lemma~\ref{NPrank-preserve} we introduce a special class of positions in $P_A$ and  a special class  of conservative $d$-pseudomoves. For a matrix $X$, we call an active column $j$ \emph{saturated} if it has exactly $r$ zeros in the normalized matrix $C(X)$ (equivalently, $|Z(j)|=r$ in the scan procedure). 
A~matrix $X$ is \emph{saturated} if all its  active columns are saturated.

To  determine the capacity it is sufficient to consider pseudomoves in a special form. 
Let $X\xrightarrow{d}Y$ be   a conservative $d$-pseudomove.
Let $m(i)$ be the smallest $j$ such that $Y_{ij}\ne X_{ij}$, i.e.\ the leftmost change in row $i$. If row $i$ is unchanged, set $m(i) = \infty$. We call $m(i)$ the \emph{first-change time} in~$i$. The leftmost change must be $1\to 0$ while the subsequent changes are arbitrary. So, $Y_{i m(i)}=0$ if $m(i)<\infty$.  
Put $m'(i)=\min\big(m(i),\ell(i)\big)$ if $m(i)<\infty$ and $m'(i)=\infty$ otherwise. Define $Y'$ by the normalization rules:
\begin{equation}\label{eq:normal-pseudomove}
 Y'_{ij}=\left\{\begin{aligned}
 X_{ij},&&&j< m'(i),\\
 0, &&& j = m'(i),\\
 \lambda(X)_j,&&&j> m'(i).
 \end{aligned}\right.
\end{equation}
Pseudomoves in the form $X\xrightarrow{d}Y'$ are called \emph{normalized}. The first-change time in row $i$ of $Y'$ is  $m'(i)\leq\ell(i)$, and $X_{im'(i)}=1$. In particular, every changed row has zeros in all inactive columns of $Y'$. Note that the first-change times  for the normalized pseudomove $X\xrightarrow{d}Y'$ equal $m'(i)$ defined for the pseudomove $X\xrightarrow{d}Y$.

\begin{proposition}\label{A-saturated}
  If $X$ is saturated then $A(X) = |E|$. 
\end{proposition}
\begin{proof}
  Normalization provides a conservative $|E|$-pseudomove (see the proof of Proposition~\ref{AgeqE}).
  Suppose for a contradiction that there exists a conservative $d$-pseudomove  $X\xrightarrow{d}Y$ with $d>|E|$.   
  Take $j=\min\{m(i):m(i)<\ell(i)\} $. The minimum is taken over a nonempty set since in this case at least one   active row is changed.   Choose a row $i$ attaining this minimum, so that $m(i) = j < \ell(i)$.
 Then $i\in I(j)$ since $j<\ell(i)$. 
 The change of the entry $(i,j)$ is $1 \to 0$. So, the column must be active,  otherwise $X_{ij}=0$. 
Since the pseudomove is conservative, $\ld(X)_j= \ld(Y)_j=1$.
 Thus, after that change an additional zero is added to zeroes in  active column $j$.
  Since the matrix is saturated, $|Z(j)|=r$.  If  $p\in Z(j)$  then  $j< \ell(p)$, hence, $m(p)\geq  j$ by the choice of  $j$. Moreover, $m(p)> j$ since $X_{pj}=0$. So the entry $(p,j)$ is unchanged. Thus, the number of zeros in active rows increases   and       column $j$ becomes inactive. The rank is changed, a~contradiction.
\end{proof}

\begin{proposition}\label{normalized-pseudomove}
 For any   conservative $d$-pseudomove   $X\xrightarrow{d}Y$, the pseudomove defined by~\eqref{eq:normal-pseudomove} is  a conservative $d$-pseudomove from $X$.
\end{proposition}
\begin{proof}
  If $m(i)=\infty$, then Eq.~\eqref{eq:normal-pseudomove} shows that row $i$ in $Y'$ coincides with row $i$ in $X$. On the other hand, if $m(i)<\infty$, then the second line of Eq.~\eqref{eq:normal-pseudomove} guarantees that row $i$ strictly decreases: either $m'(i)=\ell(i)$ and $X_{i\ell(i)}=1$ by definition of the normalized matrix, while $Y'_{i\ell(i)}=0$, or $m'(i)=m(i)$ and $X_{im(i)}=1$, $Y'_{im(i)}=0$ by definition of $m(i)$. In both cases the preceding entries are unchanged. Thus, exactly $d$ rows are strictly decreased.

  We denote by $Z(j)$, $I(j)$ the sets generated by the scan procedure on matrix $Y$ and by $Z'(j)$, $I'(j)$ the corresponding sets for matrix $Y'$. By induction on $j$, we prove that
  \begin{enumerate}
  \item $Z'(j)\supseteq Z(j)$ for an inactive column $j$ in $Y$,
  \item $Z'(j)\subseteq Z(j)$ for an active column $j$ in $Y$,
  \item $I'(j)\supseteq I(j)$ in both cases,
  \item\label{prop-norm-last} if $i\in I'(j)\sm I(j)$, then $m'(i)\leq j$.
  \end{enumerate}
The first two inclusions imply that corresponding columns receive the same rank bit, so the rank is preserved.

  Initially, the active sets are the same and claim~\ref{prop-norm-last} is vacuous. For the induction step, suppose that the claims hold for all columns preceding column $j$. Since $X\xrightarrow{d}Y$ is conservative, the rank words of $X$ and $Y$ coincide. If $m'(i)>j$, then $Y'_{ij}=X_{ij}=Y_{ij}$. If $m'(i)<j$, then $Y'_{ij}=\ld(X)_j$, and if $m'(i)=j$, then $Y'_{ij}=0$ by Eq.~\eqref{eq:normal-pseudomove}.

  Let $j$ be an inactive column in $Y$. Every row in $Z(j)$ belongs to $I'(j-1)$ by the induction hypothesis, and the preceding observations give $Y'_{ij}=0$ in each such row. Thus $Z'(j)\supseteq Z(j)$, so column $j$ is inactive in $Y'$ as well, and $I'(j)=Z'(j)\supseteq Z(j)=I(j)$. To verify claim~\ref{prop-norm-last}, take $i\in I'(j)\sm I(j)$. If $i\notin I(j-1)$, then $i\in I'(j-1)\sm I(j-1)$, so $m'(i)\leq j-1$ by the induction hypothesis. Otherwise, row $i$ is excluded from $I(j)$ at this column, so $Y_{ij}=1$ while $Y'_{ij}=0$. This excludes $m'(i)>j$ and gives $m'(i)\leq j$.

  Let $j$ be an active column in $Y$. By the induction hypothesis, $I'(j-1)\supseteq I(j-1)$, and $i\in I'(j-1)\sm I(j-1)$ implies $m'(i)<j$, hence $Y'_{ij}=1$. Thus every row in $Z'(j)$ belongs to $I(j-1)$. In such a row, either $m'(i)>j$, so $Y_{ij}=Y'_{ij}=0$, or $m'(i)=j$. In the latter case, $j$ is active in $X$, so it cannot be the marking column $\ell(i)$. Therefore $m'(i)=m(i)$ and $Y_{ij}=0$ as well. Consequently, $Z'(j)\subseteq Z(j)$ and $|Z'(j)|\leq |Z(j)|\leq r$. Column $j$ is active in $Y'$, and $I'(j)=I'(j-1)\supseteq I(j-1)=I(j)$. Claim~\ref{prop-norm-last} follows from the induction hypothesis.
\end{proof}

\begin{proof}[Proof of Lemma~\ref{NPrank-preserve}]
  Let $X\in N_A$. Thus, $A(X)\geq k$.   Let  $\lambda=\lambda(X)$, $I$, and $E$  be the rank word, the active and
marked row sets of $X$, respectively.

The first case: $|E|\geq k$. Select a set $S$ of $k$ marked
rows and apply normalization~\eqref{eq:normal} to these rows,
leaving all other rows unchanged. This preserves the rank:
a normalized row $i$ remains active after its marking column $\ell(i)$
and subsequently follows $\lambda$, contributing a one in
each active column and a zero in each inactive column.
The normalized matrix of the resulting position is $C(X)$,
and its marked rows are exactly $E\setminus S$.

Now saturate each active column by changes $1\to0$ in the
selected rows. To justify this step, note that $r=n-2k\leq k$ since we assume that $n\leq 3k$. Let $z\leq r$ be the number of zeros in the corresponding column of $C(X)$.
At least $k-z$ of its entries in the selected rows are ones.
Since $  k-z\geq r-z$, 
we can change $r-z$ of these entries to zero. Do this
independently in every active column, and denote the
resulting position by $Y$.

These additional changes preserve the rank. Indeed, a
column-by-column comparison with the scan of $X$ shows
that the active row set after column $j$ is $I(j)\cup S$.
At an inactive column the selected rows are zero, and the
unselected rows follow the original scan. At an active
column the zeros among the active rows are exactly the
$r$ zeros in the modified normalized column. Hence the
rank word remains $\lambda$, and the marked rows of $Y$
are exactly $E\setminus S$. Every selected row was already
strictly decreased by normalization, and the additional
changes only decrease it further. Thus $X\to Y$ is a legal
move. The matrix $Y$ is saturated, so
Propositions~\ref{E-bound} and~\ref{A-saturated} give $A(Y ) = |E|-k <k$ since $|E|<2k$ by Proposition~\ref{E-bound}.
Therefore $Y\in P_A$ and $b(Y)=b(X)$.

The second case: $|E|<k$. By
Proposition~\ref{normalized-pseudomove}, choose a normalized
conservative pseudomove $  X\xrightarrow{A(X)}Y'$,
and let $m(i)$ be its  first-change times. Thus
$m(i)\leq\ell(i)$ whenever $m(i)<\infty$. Every marked row
is changed: otherwise, normalizing an unchanged marked row
would preserve the rank and increase the number of
decreased rows, contradicting maximality. Consequently,
every row of $Y'$ is zero in each inactive column of $X$.
Since $b(Y')=b(X)$, every active column of $Y'$ contains at
most $r$ zeros; otherwise, the OR of $r+1$ rows having zero
there would be strictly smaller than $\lambda$.

Select all marked rows and $k-|E|$ active rows changed by
the pseudomove, taking the latter in increasing order of
$m(i)$, with arbitrary resolution of ties. Let $S$ be the
resulting set of $k$ rows, and put
\[
  m=\max\{m(i):i\in S\cap I\}.
\]
This maximum is finite because $|E|<k$ and $A(X)\geq k$.
Every unselected active row has first-change time at
least $m$, possibly equal to infinity.

Define $Y''$ as follows:
\begin{itemize}
\item if $i\notin S$, leave row $i$ unchanged;
\item if $i\in S$ and either $i\in I$ or $m(i)\leq m$,
      copy row $i$ of $Y'$;
\item if $i\in  E\subseteq S$ and $m(i)>m$,  apply normalization to  $i$.
\end{itemize}
Each selected row strictly decreases, and all unselected
rows remain unchanged. Hence $X\to Y''$ is a legal move.
The  first-change times of this move are
\[
  \widehat m(i)=\left\{
  \begin{aligned}
    &\infty,
       && i\notin S,\\
    &m(i),
       && i\in S\cap I
         \text{ or }(i\in  E\ \text{and}\ m(i)\leq m),\\
    &\ell(i),
       && i\in  E,\ m(i)>m.
  \end{aligned}\right.
\]

We next verify that this move preserves the rank.
Every row of $Y''$ is zero in the inactive columns of $X$:
all marked rows were selected, and all unselected rows
belong to $I$. For an active column $j<m$, the columns of
$Y''$ and $Y'$ coincide. At $j=m$, they can differ only
because an unselected active row with $m(i)=m$ has not
been decreased. Such an omission replaces a zero of $Y'$
by a one, so the number of zeros cannot increase. Thus,
for every active column $j\leq m$, the number of zeros
in $Y''$ is at most that in $Y'$, and hence at most $r$.

For an active column $j>m$, every row copied from $Y'$
has entry $\lambda_j=1$, because its first change
occurred at or before $m$. Each remaining row of $Y''$
coincides with the corresponding row of $C(X)$: it is
either an unselected active row or a selected marked
row to which normalization was applied. Therefore the
number of zeros in such a column is at most that in
$C(X)$, again at most $r$. It follows that the OR of
every $r+1$ rows of $Y''$ is $\lambda$, and consequently
$  b(Y'')=b(X)$.

Finally, saturate every active column of $Y''$ by changing
ones to zeros only in the selected rows. If a column has
$z\leq r$ zeros, at least $k-z\geq r-z$ selected entries
are ones, so the required $r-z$ changes are possible.
The zeros at the first-change entries
$(i,\widehat m(i))$ are not altered. Additional changes
$1\to0$ cannot destroy the strict decrease of a selected
row, and no unselected row is changed.

Denote the resulting position by $Y$. Its inactive
columns are zero and every active column has exactly
$r$ zeros. Thus every $r+1$ rows have OR $\lambda$,
so $b(Y)=b(X)$. Moreover, $Y$ is normalized and saturated,
with no marked rows. Proposition~\ref{A-saturated} gives $ A(Y)=0<k$.
Therefore $Y\in P_A$, and the legal move $X\to Y$
preserves the rank, as required. 
\end{proof}


\section{Complexity issues}\label{sec:scan-final}

In this section we discuss the computational complexity of  
basic algorithmic problems related to Exact $\NIM(n,k)$ with $2k\leq n \leq 3k$.
An input to such a problem is a position $x$ in the game represented by the Bouton matrix $B(x)$.
The problems are: decide whether $x\in\Pset$; compute the value of the remoteness function $R(x)$; compute an optimal move.

We consider two cases: $k$ is fixed, i.e., we are solving a particular game; 
and $k$ is part of the input, i.e., we are going to solve all games  
$\NIM(n,k)$ with $2k \leq n \leq 3k$, by the same algorithm;  
in other words, $n$ and $k$ are part of the input.

\subsection{Computing capacity by a greedy algorithm}
\label{ssec:greedy}

The capacity of a position is defined as the value of an optimization problem. We show that this problem has a rather simple greedy solution. Due to Proposition~\ref{normalized-pseudomove} it is sufficient to find optimal values $m(i)$ since the normalized pseudomove is determined by Eq.~\eqref{eq:normal-pseudomove}.

In this solution we use \emph{deficits}. In the definition we use sets from the scan procedure~\eqref{eq:word-scan} defining active and marked rows.  For a row $i$ of a matrix $X$, the deficit $d(i)$ is defined by the rules
\begin{equation}\label{eq:deficit}
  d(i)_j = \left\{
  \begin{aligned}
    1,&&&\text{if}\ i\in Z(j)\ \text{and}\ \ld_j =1,\\
    0,&&&\text{otherwise.}
  \end{aligned}\right.
\end{equation}
Eq.~\eqref{eq:word-scan} immediately implies that the deficit matrix $(d(i)_j)$ has at most $r$ ones in each column and exactly $r$ ones in a saturated column. 
Deficits and their suffixes are compared by their binary values (equivalently, in lexicographic order).

Now we present the enhanced scan procedure generalizing the scan~\eqref{eq:word-scan}. It maintains the same sets $I(j)$, $Z(j)$, $E(j)$ as~\eqref{eq:word-scan} does and, additionally, the values $m(i)(j)$. 
Initially, $I(0)=\{1,\ldots,n\}$, $E(0)=\es$, and $m(i)(0)=\infty$ for all $1\leq i\leq n$.

 At column $j$, first repeat the actions of the scan~\eqref{eq:word-scan}. Set $Z(j)=\{i\in I(j-1):X_{ij}=0\}$ and, for $1\leq j\leq L$, set
\begin{equation*}
\left\{ \begin{aligned}
 \lambda_j&=0, &&I(j) = Z(j), && E(j) = \{1,2,\ldots, n\}\sm Z(j),&&\text{if } |Z(j)|\geq r+1,\\
 \lambda_j&=1,&&I(j) = I(j-1), && E(j) = E(j-1),&&\text{if } |Z(j)|\leq r.
 \end{aligned}\right.
\end{equation*}
After that, let $M(j) = \{i: m(i)(j-1)<\infty\}$ be the set of selected rows.

If $\ld_j=0$ then set  $m(i)(j) = j$ for each $i\in E(j)\sm E(j-1)$ with $i\notin M(j)$ and $m(i)(j)= m(i)(j-1)$ otherwise.  


If $\ld_j =1$ then compute $m(i)(j)$ as follows.
Let  $F(j) = \{i\in I(j):X_{ij}=1\}\setminus M(j) $. If  $|Z(j)\setminus M(j)|=r$ then $m(i)(j)= m(i)(j-1)$ for all $i$. Otherwise select  $p = \min (|F(j)|, r-|Z(j)\setminus M(j)|)$ rows from the set $F(j)$. The choice is made in decreasing order of the values of deficit suffixes that start from $j+1$.
For equal suffixes, the decreasing order of $\ell(i)$ is applied. In particular, active rows with $\ell(i)=\infty$ are selected before marked ones with $\ell(i)<\infty$. Remaining ties are broken arbitrarily.
Set  $m(i)(j) = j$ for a selected $i$ and $m(i)(j)= m(i)(j-1)$ otherwise.  

Finally, we set $I = I(L)$, $E= E(L)$, $m^*(i) = m(i)(L)$,
$M^*=M(L+1) = \{i:m^*(i)<\infty\}$ and $d = |M^*|$.

The following example illustrates the definition.


Take $k=3$, $r=2$, $n=8$, and $ x=(29,31,61,46,61,51,61,59)$.
Its Bouton matrix, the normalization and the column of $\ell(i)$ are
\[
 X=\begin{pmatrix}
 0&1&1&1&0&1\\
 0&1&1&1&1&1\\
 1&1&1&1&0&1\\
 1&0&1&1&1&0\\
 1&1&1&1&0&1\\
 1&1&0&0&1&1\\
 1&1&1&1&0&1\\
 1&1&1&0&1&1
 \end{pmatrix},\qquad
 C=\begin{pmatrix}
 0&1&1&1&0&1\\
 0&1&1&1&0&1\\
 1&1&1&1&0&1\\
 1&0&1&1&0&1\\
 1&1&1&1&0&1\\
 1&1&0&0&0&1\\
 1&1&1&1&0&1\\
 1&1&1&0&0&1
 \end{pmatrix}
,\qquad
\ell = 
\begin{pmatrix}
  \infty\\ 5 \\ \infty\\ 5\\ \infty\\ 5\\ \infty\\ 5
\end{pmatrix} .
\]
The ordinary scan gives
\[
 \lambda(x)=111101,\qquad b(x)=61,\qquad
 I(L)=\{1,3,5,7\},\qquad E(L)=\{2,4,6,8\}.
\]
The deficit words are
\[
 \begin{array}{c|cccccccc}
 i&1&2&3&4&5&6&7&8\\\hline
 d(i)&100000&100000&000000&010000&000000&001100&000000&000100
 \end{array}.
\]
To distinguish the two stages of a column, write
$M_j^- =M(j)$ and $M_j^+=\{i:m(i)(j)<\infty\}$; thus $M_{j+1}^-=M_j^+$.
Put $q_j=|Z(j)\setminus M_j^-|$. In active columns the available number
of new zeros is $r-q_j$. Remaining ties are resolved by increasing row
index. The complete trace is as follows; a dash denotes a quantity not
used in that branch.
\begin{center}
\small
\setlength{\tabcolsep}{4pt}
\begin{tabular}{@{}ccclclcl@{}}
\toprule
$j$&$\lambda_j$&$Z(j)$&$M_j^-$&$q_j$&$F(j)$&New rows&$M_j^+$\\
\midrule
1&1&$\{1,2\}$&$\varnothing$&2&$\{3,4,5,6,7,8\}$&$\varnothing$&$\varnothing$\\
2&1&$\{4\}$&$\varnothing$&1&$\{1,2,3,5,6,7,8\}$&$\{6\}$&$\{6\}$\\
3&1&$\{6\}$&$\{6\}$&0&$\{1,2,3,4,5,7,8\}$&$\{1,8\}$&$\{1,6,8\}$\\
4&1&$\{6,8\}$&$\{1,6,8\}$&0&$\{2,3,4,5,7\}$&$\{3,5\}$&$\{1,3,5,6,8\}$\\
5&0&$\{1,3,5,7\}$&$\{1,3,5,6,8\}$&--&--&$\{2,4\}$&$\{1,2,3,4,5,6,8\}$\\
6&1&$\varnothing$&$\{1,2,3,4,5,6,8\}$&0&$\{7\}$&$\{7\}$&$\{1,\ldots,8\}$\\
\bottomrule
\end{tabular}
\end{center}

Column 1 is the blocked case $q_j=r$. Column 2 is a genuine partial
selection: the suffix $d(6)[3,6]=1100$ is larger than every competing
suffix, including $d(8)[3,6]=0100$. In column 3, row 8 is selected first
because its remaining deficit suffix is $100$. The other candidate
suffixes are zero, so final active rows, whose marking time is infinite,
precede rows 2 and 4, whose marking time is 5; the row-index tie chooses
row 1. In column 4 all candidate deficit suffixes are zero. The marking-time
tie-break puts rows $3,5,7$ before rows $2,4$, and the remaining tie-break chooses
rows 3 and 5. In column 5, rows $2,4,6,8$ become marked: rows 6 and 8 were
already selected, while rows 2 and 4 are selected now. Finally, column 6
realizes the nonempty strict-underfill case
\[
 0<|F(6)|=1<r-|Z(6)\setminus M(6)|=2.
\]
The single candidate is selected even though a second zero could be accommodated.

Finally, $ m^*=(3,5,4,5,4,2,6,3)$ and the result of the corresponding normalized pseudomove is
$y^*=(21,29,57,45,57,45,60,53)$.
The bitwise OR of any three rows of $B(y^*)$ is $111101$. Hence
$b(y^*)=b(x)=61$ and $A(x)=8$, since all eight piles decrease.
This is a conservative eight-pseudomove, not a legal move
of $\NIM(8,3)$.

The marking-time tie-break is essential. In column 4 all
candidate deficit suffixes are zero, but rows $3,5,7$ have
$\ell(i)=\infty$, whereas rows $2,4$ have $\ell(i)=5$.
Selecting rows $2,4$ instead of $3,5$ would produce no new
selection in column 5 and leave three candidates for only
two places in column 6. The resulting pseudomove would
therefore decrease only seven piles.

\begin{lemma}\label{lem:capacity-count}
  $A(X) = d$.  
\end{lemma}
\begin{proof}
A normalized pseudomove is determined by the values $m(i)$ of its leftmost changes. So, we indicate such a move by these values.

By construction, applying the rules~\eqref{eq:normal-pseudomove} to $m^*(i)$ gives a normalized conservative $d$-pseudomove, since the status of each column remains unchanged. Thus, $A(X)\geq d$. It remains to prove the optimality of the pseudomove determined by $m^*(i)$.

For optimality, we prove by induction on $j$ that there exists an optimal normalized pseudomove $\tdm(i)$ consistent with $m^*(i)$ up to column $j$. This means that
$\{i:\tdm(i)\leq j\}=\{i:m^*(i)\leq j\}$ and $\tdm(i)=m^*(i)$ on this set. At $j=L$, the consistency condition implies $m^*=\tdm$, so the pseudomove $m^*$ is optimal.

The base case is $j=0$. In this case the consistency condition holds trivially, and an optimal normalized pseudomove exists by Proposition~\ref{normalized-pseudomove}.

For the induction step, assume that $m^*(i)$ is consistent with some optimal normalized pseudomove determined by $\tdm(i)$ up to column $j-1$.

Suppose that $\ld_j=0$. Let $i$ be a row such that $m(i)(j)=j$. Then, by construction of the enhanced scan, $i\in E(j)\sm\big(E(j-1)\cup M(j)\big)$. From the induction hypothesis, $\tdm(i)\geq m^*(i)=j$. Set $m'(i)=j$ for $i\in E(j)\sm\big(E(j-1)\cup M(j)\big)$ and $m'(i)=\tdm(i)$ otherwise. The function $m'$ determines a normalized conservative pseudomove, since the normalization rule preserves the status of every column. The number of decreased rows does not decrease, so optimality of $\tdm$ implies that $m'$ is also optimal and is consistent with $m^*$ up to column $j$.

Suppose that $\ld_j=1$ and $|Z(j)\setminus M(j)|=r$. In this case there is no room for a change $1\to0$ in an unselected row, since any such change makes column $j$ inactive and changes the rank. Thus, $\tdm(i)\ne j$ for every unselected row $i$, and the consistency condition holds at column $j$.

Suppose that $\ld_j=1$, $|Z(j)\setminus M(j)|<r$, and $|F(j)|\leq r-|Z(j)\setminus M(j)|$. In this case each row $i\in F(j)$ receives the value $m(i)(j)=j$ and becomes selected. From the induction hypothesis, $\tdm(i)\geq j$ for $i\in F(j)$. Let $S=\{i\in F(j):\tdm(i)>j\}$, and set $m'(i)=j$ for $i\in S$ and $m'(i)=\tdm(i)$ otherwise. This determines a conservative pseudomove: column $j$ remains active, and after the new first change the normalization rule puts ones in active columns and zeros in inactive columns. The number of decreased rows does not decrease. Therefore, optimality of $\tdm$ implies that $m'$ is also optimal and consistent with $m^*$ up to column $j$.

Now consider the case $\ld_j=1$, $|Z(j)\setminus M(j)|<r$, and $|F(j)|>r-|Z(j)\setminus M(j)|=p$. In this case only $p$ rows become selected while scanning column $j$. Let $\tS=\{i\in F(j):\tdm(i)=j\}$. This set can differ from the corresponding subset $S^*=\{i:m^*(i)=j\}$.
Take an optimal normalized $\tdm$, consistent with $m^*$ up to column $j-1$, for which the cardinality of the symmetric difference $\tS\sdif S^*$ is smallest. To complete the induction step, we prove that $|\tS\sdif S^*|=0$, that is, $\tS=S^*$.

Suppose for a contradiction that $|\tS\sdif S^*|>0$. The zero count in column $j$ gives $|\tS|\leq p$. If $|\tS|<p$, it is possible to decrease $\tdm(i)$ to $j$ for some $i\in S^*\sm\tS$, preserving the conditions for a conservative pseudomove, as in the preceding case. The resulting pseudomove is also optimal, while $|\tS\sdif S^*|$ decreases, a contradiction. Thus, $|\tS|=p$.

Consequently, there are $a,b\in F(j)$ such that $\tdm(a)=j$, $m^*(a)>j$, $\tdm(b)>j$, and $m^*(b)=j$.
Define $m'(i)=\tdm(i)$ for $i\ne a,b$, and set $m'(b)=j$. The selection rule guarantees that
$\val(d(a)[j+1,L])\leq\val(d(b)[j+1,L])$.
If these suffixes differ, let $q$ be their first differing column; then $d(a)_q=0$ and $d(b)_q=1$, so column $q$ is active. If the suffixes coincide, put $q=\infty$. In both cases set
\[
  m'(a)=\tau=\min\{\tdm(b),\ell(a),q\}.
\]
All three quantities in this minimum are greater than $j$. If $q<\infty$, then $d(b)_q=1$ implies $\ell(b)>q\geq\tau$. If $q=\infty$, the marking-time tie break gives $\ell(b)\geq\ell(a)\geq\tau$.

The number of decreased rows does not decrease under this reassignment. Indeed, if $\tau=\infty$, then $\tdm(b)=\infty$, and the decrease of row $a$ is replaced by a decrease of row $b$. If $\tau<\infty$, both rows are decreased by $m'$, whereas among these two rows $\tdm$ decreases either one or both. All other rows have the same first-change times.

We next verify that the new first change in row $a$ is legal when $\tau<\infty$. If $\tau=\ell(a)$, then $X_{a\tau}=1$ by the definition of the marking column. If $\tau<\ell(a)$ and $\tau=q$, then column $q$ is active and $d(a)_q=0$, so again $X_{a\tau}=1$. In the remaining case, $\tau=\tdm(b)<\min\{\ell(a),q\}$. The inequalities above imply $\tau<\ell(b)$, so the original first change in row $b$ occurs in an active column. Since the deficit bits of the two rows agree before $q$, we have $X_{a\tau}=X_{b\tau}=1$. Thus, the prescribed change in row $a$ is $1\to0$. The change in row $b$ at column $j$ is legal because $b\in F(j)$. All finite first-change times remain at most the corresponding marking times, so $m'$ defines a normalized pseudomove.

It remains to verify that this pseudomove is conservative. Columns preceding $j$ are unchanged, and at column $j$ the new zero in row $b$ replaces the old zero in row $a$. For $j<t<\tau$, both rows remain active, and $d(a)_t=d(b)_t$. In an active column $t$, delaying the first change in row $a$ introduces exactly as many zeros as are removed by advancing the first change in row $b$. In an inactive column $t$, both original entries are zero, since $t<\ell(a),\ell(b)$, and the normalized entries are also zero.

If $\tau<\infty$, consider column $\tau$. When $\tau=\ell(a)$, this column is inactive and the changed rows contribute zeros, so it remains inactive. When $\tau<\ell(a)$ and $\tau=q$, the original entry $X_{bq}=0$ implies $\tdm(b)>q$; hence the zero introduced in row $a$ replaces the old zero in row $b$. In the remaining case, $\tau=\tdm(b)<\min\{\ell(a),q\}$, the zero formerly introduced by the first change in row $b$ is replaced by the first-change zero in row $a$.
After column $\tau$, both reassigned rows follow the rank word. Thus, compared with the pseudomove $\tdm$, the zero count cannot increase in an active column or decrease in an inactive column. Any additional active row contributes only rank-word entries. A column-by-column comparison therefore shows that every column retains its status and the rank is preserved. If $\tau=\infty$, the preceding comparison for $j<t<\tau$ covers all remaining columns.

We have obtained a conservative pseudomove whose number of decreased rows is at least that of the optimal pseudomove $\tdm$. By optimality, these numbers are equal, so $m'$ is also optimal. Since $m'(b)=j$ and $m'(a)>j$, the symmetric difference with $S^*$ decreases, a contradiction. Therefore $\tS=S^*$, and the induction step is complete.
\end{proof}

\begin{remark}\label{72-nonmonotone}
We define the capacity by means of the rank function. 
The following monotonicity property of N-positions is essential:  
Given an N-position, add to its Bouton matrix 
arbitrary columns on the right, 
then the resulting extended matrix also defines an N-position.      

This property holds when 
$2k \leq n \leq 3k$ since the capacity is not decreased after adding less significant columns by Lemma~\ref{lem:capacity-count}: note that the enhanced scan  never unselects rows and   any choices in the enhanced scan generate a normalized conservative pseudomove.
But it  does not hold for all possible $n$, $k$. 
For example, in $\NIM(7,2)$ we have: 
$(0,0,0,0,1,1,1)$ is an N-position, while   
$(1,1,1,1,2,2,2)$ is a P-position.
\end{remark}

\subsection[Parameters n and k as part of the input]{Parameters \boldmath $n$ and $k$ as part of the input}

Now we assume that $n$ and $k$ are part of the input and $2k\leq n\leq 3k$.

\begin{theorem}\label{thm:nk-input}
   Recognizing P-positions, computing the remoteness function,  
   and finding an optimal move from a nonterminal position can be done in polynomial time.
\end{theorem}
\begin{proof}
  To implement the enhanced scan procedure one needs to maintain the sets used in the procedure. This can be done easily in polynomial time and the most time-consuming operation is the ordering of $F(j)$ according to the values of deficit suffixes.
This proves the first two claims. The optimal move from a nonterminal P-position in Lemma~\ref{lem:P-Smith} can be constructed in polynomial time provided the rank and the set of active rows have already been computed.

  An optimal move from an N-position can be computed following the proof of Lemma~\ref{NPrank-preserve}. For an optimal pseudomove needed in the proof, one can take the pseudomove generated by the enhanced scan. 
\end{proof}


\subsection{A fixed game}

Let $k$ be a fixed integer, $n=2k+r$, $0\leq r\leq k$. The computational complexity of basic algorithmic problems related to $\NIM(n,k)$ describes the intrinsic complexity of the game. So, it is interesting to give tighter bounds on computational complexity in this case.

For finite-state automata and finite-state transducers we adopt definitions from Berstel's book~\cite{Berstel}.

The enhanced scan procedure described in Section~\ref{ssec:greedy} can
be converted to a nondeterministic finite automaton (NFA) $\A$  reading columns of the Bouton matrix. So the alphabet of the automaton is  $\{0,1\}^n$. It consists of $2^{n}$ symbols.  The automaton $\A$ accepts exactly Bouton matrices $B(x)$ such that $x\in\Nset$.

States of the automaton $\A$ carry the sets used in the enhanced scan procedure, namely $I$ and $M$.
There are $4^n$ pairs $(I,M)$.

The initial state is  $I=\{1,\ldots,n\}$, $ M =\es$.
A state $(I,M)$ is accepting if and only if $|M|\geq k$. 
Transitions are defined as follows. Let $c\in\{0,1\}^n$ be a column. From $c$ and the current state $(I, M)$, the sets $E=\{1,\dots, n\}\sm I$ (the marked-row set),  $Z=\{i: i\in I\ \text{and}\ c_i=0 \}$ (the zero set),  $F = \{i\in I: c_{i}=1\}\setminus M $ (the candidates for selection) are restored.
The updates depend on $|Z|$ and follow the rules of the enhanced scan with the following modifications. The automaton $\A$ does not maintain $m(i)(j)$. If an assignment $m(i)(j)= j$ is performed in the enhanced scan, $\A$ simply adds $i$ to the set $M$. If a choice of $p$ rows from the set $F(j)$ is performed, $\A$ nondeterministically guesses the selected subset.  

\begin{remark}
Some states of $\A$ are not reachable from the initial state.
If a row is marked then it is selected at some step of the enhanced scan. So $\{1,\dots,n\}\sm I\subseteq M$ in every reachable state. Thus $\A$ has at most $3^n$  states that are reachable from the initial state.
\end{remark}

\begin{proposition}\label{prop:scan}
  Let $L(\A)\subseteq \big(\{0,1\}^n\big)^*$ be the language recognized by the automaton $\A$. Identify words in $\big(\{0,1\}^n\big)^*$ with $n$-row matrices.  Then $L(\A) = \Nset$.
\end{proposition}
\begin{proof}
  Any run of $\A$ on an input  $c = c_1c_2\dots c_L\in  \big(\{0,1\}^n\big)^*$ follows the enhanced scan except for the choice of $p$ rows from $F(j)$. The states along the run  $(I_j, M_j)$, where $0\leq j\leq L$, are  $I_j = I(j)$, $M_j= M(j+1)$ for $1\leq j\leq L$. Thus, $M^* = M_L$. For $M_j$, note that an assignment $m(i)(j)= j$ is applied in the enhanced scan only to the rows  such that $m(i)(j-1)=\infty$  and all these rows are added to $M(j+1)$.

  There exists a run making correct guesses at each step requiring a guess. The run finishes with $|M^*|= A(X)$. So, it accepts any input from $\Nset$.

  For the inputs from $\Pset$, note that any choices in the enhanced scan generate a normalized conservative pseudomove. So, $|M^*|<k$ for any run of $\A$ on an input from  $\Pset$.
\end{proof}

\begin{remark}
  Proposition~\ref{prop:scan} implies that the predicate $P(x_1,\dots,  x_n)$ which is true if $(x_1,\dots,x_n)$ is a P-position in Exact $\NIM(n,k)$, $2k\leq n\leq 3k$, can be defined by a first-order formula in B\"uchi arithmetic $\BA_2$. It is a structure $\big(\NN, x=y, x+y, V_2(x)\big)$, where $=, +$ have the usual meaning and $V_2(x)$ is  the largest power of $2$ dividing $x$ if $x \ne 0$; $V_2(0)=1$. For details on B\"uchi arithmetic and its relation to automata reading digits of integers, see~\cite{BHMV94}.
\end{remark}

Recall that for any NFA there exists an equivalent deterministic finite automaton (DFA) and  
regular languages are closed under taking the complement. So, there exist DFAs $\A_\Pset$ recognizing the set $\Pset$ and $\A_\Nset$   recognizing the set $\Nset$.  Note also, that restricting $\A$ to the ordinary scan gives a~DFA $\A_I$ that maintains the active row sets $I(j)$ while reading the input (all states of  $\A_I$ are accepting).
2
  The DFAs $\A_\Nset$ and $\A_I$ are easily converted to  nondeterministic finite-state transducers (NFST) $\TR$ and $\TS$ that output the rank word and the remoteness function. The transducer $\TR$
  runs  $\A_I$  and outputs $1$ while reading a   column such that $|Z|\leq r$ and 0 otherwise. The transducer $\TS$ runs $\A_I$ and $\A_\Nset$ in parallel,
nondeterministically guesses the end of the input, and adds to the output the bit   $\uu[x\in \Nset]$ determined by a state of $\A_\Nset$.

  It is easy to see that  $\A_\Pset$, $\TR$ and $\TS$  can be simulated in linear time. For linear-time algorithms we adopt the word RAM model (registers have $O(\log(N+1))$ bits, where $N=nL$ is the input length); see~\cite{donkey-book}. 

  Thus, for a fixed game  Theorem~\ref{thm:nk-input} can be improved: recognizing P-positions and computing the remoteness function can be done in linear time. For optimal moves the situation is more complicated. But, again, we are able to improve Theorem~\ref{thm:nk-input} and present an algorithm computing optimal moves in linear time. 

\begin{proposition}\label{Smith-optimal:NFST}
  There exists an NFST $T$ such that on any input $x\in \big(\{0,1\}^n\big)^*$ the transducer produces all optimal moves. More exactly, let $T(x)$ be the set of all possible outputs of $T$ along accepting runs on an input $x$. Then  $T(x) = \{y\colon x\to y\ \text{is optimal}\}$.
\end{proposition}
\begin{proof}
  The required transducer $T$ guesses columns of $y$ nondeterministically after reading a column from the input and checks the correctness of its guesses. If some check fails, $T$ goes to an absorbing rejecting state (thus it produces no accepted output on this run).
  All checks are done in parallel.
  
  Verifying legality of the guessed move is easy. The transducer $T$ keeps  one comparison status per row: equal so far, already smaller, or invalid (an incorrect guess). At the end exactly $k$ rows must be smaller. Note that for a terminal input this check necessarily fails. So, $T(x)=\es$ if $x$ is terminal. 

The transducer $T$ runs in parallel two copies of the automaton $\A_\Pset$  on $x$ and $y$ and two copies of $\TR$. For the first two copies $T$ guesses the result of operation of $\A_\Pset$ in advance (and checks the correctness of these guesses at the end of input). If guessed outcomes are equal then $T$  rejects. For the remaining cases, $T$ checks equalities
\[
 \begin{array}{ll}
\text{if}\ x\in\Nset,\ y\in\Pset,&\text{then}\ b(y)=b(x),\\
\text{if}\  x\in\Pset,\ y\in\Nset,&\text{then}\  b(y)=b(x)-1.
 \end{array}
\]
Equality is verified by component-wise comparison of bits produced in the computation of the rank word.  Binary decrement  by 1 can also be verified by a finite-state machine:
the rank words first agree, then have a position such that $\ld(x)_j=1$, $\ld(y)_j=0$, and the opposite equalities   $\ld(x)_p=0$, $\ld(y)_p=1$ for all positions $p>j$.
\end{proof}

The transducer in Proposition~\ref{Smith-optimal:NFST} is intrinsically nondeterministic. So, it is unclear how to compute at least one optimal move in deterministic linear time. The uniformization theorem (see~\cite{Ei74}) helps. It claims that for each NFST $T$ there exists an NFST $T'$ such that $T'$ computes a~function consistent with outputs of $T$. More exactly, $|T'(x)|\leq1$ and, if $T(x)\ne\es$, then $T'(x)=\{y\}$ for some $y\in T(x)$. Nevertheless,  $T'$ is nondeterministic in general (but all accepting paths on a specific input generate the same output).

\begin{theorem}\label{thm:Smith-move-linear-time-linear-space}
  Let $k=O(1)$, $n=2k+r$, $0\leq r\leq k$. There exists an algorithm running in linear time such that,
  given a position $x$, the algorithm either finds an optimal move $x\to y$ or reports that $x$ is terminal.
\end{theorem}
\begin{proof}
  Let $T$ be the transducer from Proposition~\ref{Smith-optimal:NFST} and $T'$ be its uniformization.
  Let $Q$ be the state set of $T'$ and $A$ be the set of  accepting states of $T'$. For the given input columns
$c_1,\ldots,c_L$, compute backward the sets $A_j$ of states from which reading  $c_{j+1}\cdots c_L$ can lead to an accepting state.
Set $A_L=A$ and use the recurrence 
\[
 A_{j-1}=\{q\in Q:\ \text{reading } c_j \ \text{leads to a state in}\  A_j\}
 \]
 to compute the sets $A_j$ backwards. 
 Linear space is needed to store all $A_j$ since each of them has size $O(1)$.

 If $q_0\notin A_0$, where $q_0$ is the initial state, the output is empty. By  the properties of $T$, this occurs exactly when the input position is terminal.

If  $q_0\in A_0$ then  the forward pass is started. Reading column $j$, it selects  a transition from the current state  to $A_j$ and writes the output of $T'$ during this transition.
\end{proof}

Regarding circuit complexity, we are able to place the problems under consideration at low levels of the $\NC$ hierarchy. Recall that  logspace-uniform $\NC^1$ is the class of Boolean functions computed by   logspace-uniform bounded-fan-in Boolean circuits of size $\poly(N)$ and depth $O(\log(N+1))$, where $N$ is the number of variables. Logspace-uniform  $\NC^1$ is contained in the class $\L$ of logarithmic space.
A related class $\FNC^1$ consists of Boolean maps that are computed by logspace-uniform multi-output circuits of size $\poly(N)$ and depth $O(\log(N+1))$.
For details, see~\cite{Vollmer}. 

\begin{theorem}\label{thm:fixed-circuits}
  Let $k=O(1)$, $n=2k+r$, $0\leq r\leq k$.  P-membership belongs to  $\NC^1$. Computing the remoteness function and an optimal move belongs to $\FNC^1$. Moreover, there are linear-size circuits computing them. 
\end{theorem}
\begin{proof}
The number of Boolean variables for these problems is $N= nL$, where $L$ is the number of columns in the input matrix. If $k=O(1)$ then  $N = \Theta(L)$. 
  
  For a fixed DFA, each input column specifies a map from a
constant-size state set to itself. Composition of these maps is
associative and can be implemented by a constant-size Boolean circuit.
Parallel prefix computation (see~\cite{LF80} and~\cite[Ch.~29]{donkey-book}) computes all prefix compositions with
$O(L)$ compositions and depth $O(\log(L+1))$. This suffices for P-membership and computation of the remoteness function.

For an optimal  move, use  Boolean relations on the state set $Q$ of the transducer $T'$ for each column.  Relation composition is again an associative operation on a constant-size object. A~reversed prefix scan
therefore computes all the suffix-feasibility sets $A_j$ in
Theorem~\ref{thm:Smith-move-linear-time-linear-space}. Once those sets are known, choose for each column $j$ a deterministic map  from $A_{j-1}$ to $A_{j}$ consistent with the transition relation of $T'$.   A second prefix scan computes all prefix compositions of these maps and outputs of $T'$  along the corresponding run.
The two scans have total size $O(L)$ and depth $O(\log(L+1))$.
A~separate terminal flag indicates that no legal move is available (the empty output of $T'$).

For fixed $k$, all finite tables are constants. Tree indices and gates in the parallel prefix computation procedures can be generated with $O(\log L)$ workspace, giving the stated uniformity.
\end{proof}

\section*{Conclusion}

We have characterized the P-positions and determined the Smith remoteness
function of Exact $\NIM(n,k)$ for $2k\leq n\leq 3k$ in terms of the rank
and capacity of the Bouton matrix. The same description provides optimal
moves from both P-positions and N-positions, while the greedy capacity
algorithm makes these results polynomial-time computable even when $n$
and $k$ are part of the input. For each fixed game, the finite-state
representation yields stronger bounds, including linear-time algorithms
and logarithmic-depth circuits of linear size. A natural direction for
further research is to determine whether a comparable description can
be obtained for $n>3k$. The example in Remark~\ref{72-nonmonotone} shows that the behavior
under the addition of less significant binary columns can change beyond
the range considered here, so an extension must account for this
difference. Another question is whether the rank and capacity, possibly
together with additional invariants, can contribute to the determination
of nonzero Sprague--Grundy values for Exact $\NIM(n,k)$ with $n>2k$.

\section*{Statement about AI usage} 

ChatGPT (GPT-6 Astra Pro) was used during the preparation of the manuscript to assist
with coding and exploratory computations. After using this tool, the authors reviewed, verified, and edited all generated outputs as needed and take full responsibility for the content of the paper.
Claude (Fable 5.0, 5.2) was used to run experiments on NIMs with fixed $(n, k)$ to identify patterns.
This model was also used to check proofs and grammar. 

\section*{Acknowledgments} 
This paper is an output of the research project (HSE-BR-2025-024) 
implemented as part of the Basic Research Program at HSE University.


\begin{thebibliography}{99}
\small
\setlength{\itemsep}{3pt}


\bibitem{Berstel} J.~Berstel, Transductions and Context-Free Languages, B.G. Teubner, Stuttgart, 1979.

\bibitem{B02}
C.L.~Bouton,
Nim, a game with a complete mathematical theory,
Ann. Math. (2) 3 (1901--1902) 35--39.
\url{https://doi.org/10.2307/1967631}.

\bibitem{BGHMM18}
E.~Boros, V.~Gurvich, N.B.~Ho, K.~Makino, P.~Mur\v{s}i\v{c},
On the Sprague--Grundy function of Exact $k$-Nim,
Discrete Appl. Math. 239 (2018) 1--14.
\url{https://doi.org/10.1016/j.dam.2017.08.007}.
First preprint: arXiv:1508.04484v1 (2015).
\url{https://arxiv.org/abs/1508.04484v1}.

\bibitem{BGHMM19}
E.~Boros, V.~Gurvich, N.B.~Ho, K.~Makino, P.~Mur\v{s}i\v{c},
Sprague--Grundy function of matroids and related hypergraphs,
Theoret. Comput. Sci. 799 (2019) 40--58.
\url{https://doi.org/10.1016/j.tcs.2019.09.041}.

\bibitem{BGHMM24}
E.~Boros, V.~Gurvich, N.B.~Ho, K.~Makino, P.~Mur\v{s}i\v{c},
Impartial games with decreasing Sprague--Grundy function
and their hypergraph compound,
Int. J. Game Theory 53 (4) (2024) 1119--1144.
\url{https://doi.org/10.1007/s00182-023-00850-7}.

\bibitem{BGMV24}
E.~Boros, V.~Gurvich, K.~Makino, M.~Vyalyi,
Computing remoteness functions of Moore, Wythoff, and Euclid's games,
Int. J. Game Theory 53 (2024) 1315--1333.
\url{https://doi.org/10.1007/s00182-024-00914-2}.
Preprint: arXiv:2311.02685v1 (2023).
\url{https://arxiv.org/abs/2311.02685}.

\bibitem{BHMV94}
V.~Bruy\`ere, G.~Hansel, Ch.~Michaux, R.~Villemaire, Logic and $p$-recognizable sets of integers, Bull. Belg. Math. Soc. Simon Stevin 1 (2) (1994) 191--238. 


\bibitem{CGKPV21}
N.~Chikin, V.~Gurvich, K.~Knop, M.~Paterson, M.~Vyalyi,
More about exact slow $k$-Nim,
Integers 21 (2021) \#G4, 1--14.

\bibitem{donkey-book}
T.H.~Cormen, C.E.~Leiserson, R.L.~Rivest,
Introduction to Algorithms, first ed.,
MIT Press, Cambridge, MA, 1990.

\bibitem{DH26}
M.~Dufour, S.~Heubach,
On the $\mathcal{P}$-positions of some infinite families of Slow $A$-Nim,
Integers 26 (2026) \#G3, 1--25.

\bibitem{Ei74}
S.~Eilenberg,
Automata, Languages, and Machines, Vol.~A,
Academic Press, New York, 1974.

\bibitem{Grundy}
P.M.~Grundy,
Mathematics and games,
Eureka 2 (1939) 6--8.

\bibitem{GHHC20}
V.~Gurvich, S.~Heubach, N.B.~Ho, N.~Chikin,
Slow $k$-Nim,
Integers 20 (2020) \#G3, 1--19.

\bibitem{GH15}
V.~Gurvich, N.B.~Ho,
Slow $k$-Nim,
arXiv:1508.05777 (2015).
\url{https://arxiv.org/abs/1508.05777}.

\bibitem{GMMN23}
V.~Gurvich, V.~Maximchuk, G.~Miheenkov, M.~Naumova,
On remoteness functions of $k$-NIM with $k+1$ piles
in normal and in mis\`{e}re versions,
Games 15 (6) (2024) 37.
\url{https://doi.org/10.3390/g15060037}.
Preprint: arXiv:2311.13511 (2023).
\url{https://arxiv.org/abs/2311.13511}.

\bibitem{GMN26}
  V.~Gurvich, V.~Maximchuk, M.~Naumova,
On remoteness functions of exact slow $\NIM$, $\NIM^1_=(4,2)$ and $\NIM^1_=(5,2)$,
Mod. Math. Methods 4 (2) (2026) 79--84. \url{https://doi.org/10.64700/mmm.104}

\bibitem{GMMV24}
V.~Gurvich, D.~Martynov, V.~Maximchuk, M.~Vyalyi,
On remoteness functions of exact slow $k$-NIM with $k+1$ piles,
Integers 24 (2024) \#G1.

\bibitem{GN24}
V.~Gurvich, M.~Naumova,
Screw discrete dynamical systems and their applications to exact slow Nim,
Discrete Appl. Math. 358 (2024) 382--394.

\bibitem{GN26}
V.~Gurvich, M.~Naumova,
The Sprague--Grundy Function of Slow Exact $k$-Nim with $k+1$ Piles, 
arXiv:2611.14546 (2026), 12 pp.
\url{https://arxiv.org/abs/2611.14546}. 

\bibitem{GPV23}
V.~Gurvich, A.~Parfenov, M.~Vyalyi,
Experimental study of the game Exact $\mathrm{NIM}(5,2)$,
arXiv:2311.18772v1 (2023).
\url{https://arxiv.org/abs/2311.18772}.

\bibitem{JenkynsMayberry}
T.A.~Jenkyns, J.P.~Mayberry,
The skeleton of an impartial game and the Nim-function of Moore's
$\mathrm{Nim}_k$,
Int. J. Game Theory 9 (1) (1980) 51--63.
\url{https://doi.org/10.1007/BF01784796}.

\bibitem{LF80}
R.E.~Ladner, M.J.~Fischer,
Parallel prefix computation,
J. ACM 27 (4) (1980) 831--838.

\bibitem{M10}
E.H.~Moore,
A generalization of the game called Nim,
Ann. Math. (2) 11 (3) (1910) 93--94.

\bibitem{S66}
C.A.B.~Smith,
Graphs and composite games,
J. Combin. Theory 1 (1966) 51--81.
\url{https://doi.org/10.1016/S0021-9800(66)80005-4}.

\bibitem{Sprague35}
R.~Sprague,
\"{U}ber mathematische Kampfspiele,
Tohoku Math. J. (1) 41 (1935) 438--444.

\bibitem{Sprague37}
R.P.~Sprague,
\"{U}ber zwei Abarten von Nim,
Tohoku Math. J. (1) 43 (1937) 351--354.

\bibitem{Vollmer} H.~Vollmer, Introduction to Circuit Complexity: A Uniform Approach, Springer-Verlag, Berlin, Heidelberg, 1999. 


\end{thebibliography}
\end{document}